\documentclass[12pt]{amsart}
\pdfoutput=1

\usepackage[margin=1in]{geometry}
\usepackage[utf8]{inputenc}
\usepackage{graphicx,cite}
\usepackage{amsmath,amssymb,amsthm}
\usepackage{color}
\usepackage{tikz}
\usepackage{pgfplots}
\pgfplotsset{compat=1.18}
\usepgfplotslibrary{fillbetween}
\allowdisplaybreaks

\makeatletter
\def\l@section{\@tocline{1}{12pt plus2pt}{0pt}{}{\bfseries}}
\def\l@subsection{\@tocline{2}{0pt}{2pc}{2pc}{}}
\makeatother

\makeatletter
\def\subsection{\@startsection{subsection}{2}{\z@}%
  {-3.25ex\@plus -1ex \@minus -.2ex}%
  {1.5ex \@plus .2ex}%
  {\normalfont\bfseries\boldmath}}
\def\subsubsection{\@startsection{subsubsection}{3}%
  \z@{.5\linespacing\@plus.7\linespacing}{-.5em}%
  {\normalfont\bfseries\boldmath}}
\renewcommand\paragraph{\@startsection{paragraph}{4}{\z@}%
  {3.25ex \@plus1ex \@minus.2ex}%
  {-1em}%
  {\normalfont\normalsize\bfseries}}
\makeatother

\theoremstyle{plain}
\newtheorem{thm}{Theorem}[section]

\newtheorem{cor}[thm]{Corollary}

\newtheorem{prop}[thm]{Proposition}

\theoremstyle{definition}
\newtheorem{defn}[thm]{Definition}

\theoremstyle{remark}
\newtheorem{rem}[thm]{Remark}

\theoremstyle{plain}

\numberwithin{equation}{section}

\newcommand{\thistheoremname}{}
\newtheorem{genericthm}[thm]{\thistheoremname}

\newtheorem*{genericthm*}{\thistheoremname}
\newenvironment{namedthm*}[1]
  {\renewcommand{\thistheoremname}{#1}%
   \begin{genericthm*}}
  {\end{genericthm*}}

\newcommand{\norm}[1]{\left\Vert#1\right\Vert}

\newcommand{\R}{\mathbb R}

\newcommand{\Z}{\mathbb Z}

\newcommand{\calD}{\mathcal D}

\newcommand{\dist}{\hbox{\rm dist}}

\newcommand{\calS}{\mathcal S}
\newcommand{\supp}{\textnormal{supp}}

\makeatletter
\newcommand{\vast}{\bBigg@{4}}
\newcommand{\Vast}{\bBigg@{5}}
\makeatother

\def\udot#1{\ifmmode
  \oalign{$#1$\crcr\hidewidth.\hidewidth}
  \else
  \oalign{#1\crcr\hidewidth.\hidewidth}
  \fi}

\def\beq{\begin{equation}}
\def\eeq{\end{equation}}

\makeatletter
\newcommand{\doublewidetilde}[1]{{%
  \mathpalette\double@widetilde{#1}%
}}
\newcommand{\double@widetilde}[2]{%
  \sbox\z@{$\m@th#1\widetilde{#2}$}%
  \ht\z@=.9\ht\z@
  \widetilde{\box\z@}%
}
\makeatother

\makeatletter
\def\@makefnmark{%
  \leavevmode
  \raise.9ex\hbox{%
    \fontsize\sf@size\z@\normalfont\tiny\@thefnmark}}
\makeatother

\usepackage[hidelinks]{hyperref}

\title[Sparse bounds for the trilinear Hilbert transform]
{Sparse domination and an extended $L^p$ range
for the trilinear Hilbert transform}

\author{Bingyang Hu}
\address{(Bingyang Hu) Department of Mathematics and Statistics\\
         Auburn University\\
         Auburn, Alabama, U.S.A, 36849}
\email{bzh0108@auburn.edu}
\date{October 28, 2026}


\begin{document}

\begin{abstract}
We establish sparse bounds for the trilinear Hilbert transform
$$
T(f_1,f_2,f_3)(x)
:=\operatorname{p.v.}\int_{\R}
f_1(x-t)f_2(x-2t)f_3(x-3t)\,\frac{dt}{t},
$$
with averaging exponents arbitrarily close to $2$ from above.
As a consequence, we extend OpenAI's
$L^3\times L^3\times L^3\to L^1$ bound (entry 086 in the OpenAI Research Catalog) to
$L^{p_1}\times L^{p_2}\times L^{p_3}\to L^r$
throughout the range
$$
2<p_1,p_2,p_3<\infty,
\qquad
\frac1r=\frac1{p_1}+\frac1{p_2}+\frac1{p_3},
\qquad
\frac23<r<2,
$$
and obtain corresponding weighted estimates for weights
in multilinear Muckenhoupt classes.
Our proof relies on the normalized local estimate
and stopping time argument developed in OpenAI's work.
\end{abstract}

\maketitle

\section{Introduction}

Consider the trilinear Hilbert transform
$$
T(f_1,f_2,f_3)(x)
:=\operatorname{p.v.}\int_{\R}
f_1(x-t)f_2(x-2t)f_3(x-3t)\,\frac{dt}{t},
$$
initially defined for $f_1, f_2, f_3 \in \calS(\R)$. OpenAI \cite{OpenAI} proved that\footnote{For a Lean~4
formalization of this result by Joris Roos, see
\url{https://github.com/roos-j/lean-trilinearht}.}
\begin{equation}\label{eq:L3}
\norm{T(f_1,f_2,f_3)}_{L^1(\R)}
\lesssim
\prod_{j=1}^3\norm{f_j}_{L^3(\R)}.
\end{equation}
The purpose of this note is to establish a sparse bound for $T$
and use it to extend the range of exponents in \eqref{eq:L3}.

Our main result is the following.

\begin{thm}\label{thm:main}
Let $2<p_1,p_2,p_3<\infty$ and set
\begin{equation}\label{eq:range}
\frac1r
=
\frac1{p_1}+\frac1{p_2}+\frac1{p_3}.
\end{equation}
If $2/3<r<2$, then
\begin{equation}\label{eq:bound}
\norm{T(f_1,f_2,f_3)}_{L^r(\R)}
\lesssim \prod_{j=1}^3\norm{f_j}_{L^{p_j}(\R)}.
\end{equation}
Consequently, $T$ has a unique bounded trilinear extension to
\[
L^{p_1}(\R)\times L^{p_2}(\R)\times L^{p_3}(\R)
\]
with values in $L^r(\R)$ for $p_1, p_2, p_3$ and $r$ satisfying the above constraints. 
\end{thm}

To prove Theorem~\ref{thm:main}, we establish a sparse bound
for the associated four-linear form
\begin{equation}\label{eq:continuous-form}
\Lambda(f_0,f_1,f_2,f_3)
:=
\int_{\R}f_0(x)T(f_1,f_2,f_3)(x)\,dx.
\end{equation}
The proof builds on the normalized local estimate in
\cite[Proposition~2.3]{OpenAI} and the stopping time argument in
\cite[Lemma~2.4]{OpenAI}.

For an interval $I\subset\R$ and $1\le q<\infty$, write
\[
\langle f\rangle_{q,I}
:=
\left(
\frac1{|I|}\int_I |f(x)|^q\,dx
\right)^{1/q}.
\]
A collection $\calS$ of intervals is called
$\eta$-sparse if there are pairwise disjoint measurable
sets $E_I\subset I$, $I\in\calS$, such that
$|E_I|\ge\eta|I|$. Define
\[
\Lambda_{q,\calS}(f_0,f_1,f_2,f_3)
:=
\sum_{I\in\calS}
|I|\prod_{j=0}^3\langle f_j\rangle_{q,I}.
\]
For simplicity, we fix the sparsity constant $\eta=3/4$ in the sequel.

\begin{thm}\label{prop:continuous}
There exists $2<q_*<3$ such that, for every
$q\in(2,q_*)$, all $f_1,f_2,f_3\in\calS(\R)$ and every bounded
compactly supported function $f_0$, we have
\begin{equation}\label{eq:continuous}
|\Lambda(f_0,f_1,f_2,f_3)|
\le
C_q\sup_{\calS}
\Lambda_{q,\calS}(f_0,f_1,f_2,f_3),
\end{equation}
where the supremum is over all $3/4$-sparse
collections of intervals in $\R$.
\end{thm}

Theorem~\ref{thm:main} follows from
Theorem~\ref{prop:continuous} by maximal function
estimates in the Banach range and the multilinear
extrapolation theorem of Li, Martell, and Ombrosi
\cite{LMO} for $r<1$.

\begin{rem}
The key point in Theorem~\ref{prop:continuous} is that
the sparse bound holds for every $q>2$ sufficiently close to $2$.
These bounds allow us to obtain the full range of exponents
stated in Theorem~\ref{thm:main} and the weighted estimates
in Corollary~\ref{cor:weighted} below.
\end{rem}

Next, we record two further consequences of our main results.

\begin{cor}\label{cor:infty}
Let $k\in\{1,2,3\}$ and suppose that
\[
p_k=\infty,
\qquad \textrm{and} \qquad
2<p_j<\infty,
\quad j\in\{1,2,3\}\setminus\{k\}.
\]
Define $r$ by
\[
\frac1r
=
\sum_{\substack{1\le j\le3\\j\ne k}}
\frac1{p_j}.
\]
If $1<r<2$, then $T$ admits a bounded trilinear
extension
\[
T:
L^{p_1}(\R)\times L^{p_2}(\R)\times L^{p_3}(\R)
\longrightarrow L^r(\R).
\]
\end{cor}

As a standard application of the sparse domination estimate in Theorem \ref{prop:continuous}, we obtain the following weighted estimates for the trilinear Hilbert transform $T$. We use the multilinear Muckenhoupt class from \cite{LMO}.
Recall that a \emph{weight} is a locally integrable function that is
positive almost everywhere, and we write $L^s(w):=L^s(\R,w(x)\,dx), \;  0<s<\infty$.

\begin{defn}
Let $\vec p=(p_1,p_2,p_3)$ satisfy
\[
2<p_j<\infty,
\qquad 1\le j\le3,
\qquad
\frac1r=\sum_{j=1}^3\frac1{p_j},
\qquad
\frac23<r<2.
\]
Set $\vec 2=(2,2,2,2)$.
For weights $\vec w=(w_1,w_2,w_3)$, define
\[
v:=\prod_{j=1}^3w_j^{r/p_j}.
\]
We say that $\vec w\in A_{\vec p,\vec 2}$ if
\begin{equation}\label{eq:weight-class}
[\vec w]_{A_{\vec p,\vec 2}}
:=
\sup_I
\left(
\frac1{|I|}\int_I v(x)^{\frac{2}{2-r}}\,dx
\right)^{\frac1r-\frac12}
\prod_{j=1}^3
\left(
\frac1{|I|}\int_I w_j(x)^{-\frac{2}{p_j-2}}\,dx
\right)^{\frac12-\frac1{p_j}}
<\infty,
\end{equation}
where the supremum is over all bounded intervals
$I\subset\R$ of positive length.
\end{defn}

\begin{cor}\label{cor:weighted}
Let $2<p_1,p_2,p_3<\infty$ and set
\[
\frac1r=\sum_{j=1}^3\frac1{p_j}.
\]
Suppose that $2/3<r<2$ and
$\vec w=(w_1,w_2,w_3)\in A_{\vec p,\vec 2}$.
With
\[
v:=\prod_{j=1}^3w_j^{r/p_j},
\]
we have
\begin{equation}\label{eq:weighted-bound}
\norm{T(f_1,f_2,f_3)}_{L^r(v)}
\le
C\prod_{j=1}^3\norm{f_j}_{L^{p_j}(w_j)}
\end{equation}
for $f_1,f_2,f_3\in C_c^\infty(\R)$.
The constant $C$ depends only on the exponents
and $[\vec w]_{A_{\vec p,\vec 2}}$.
Consequently, $T$ has a unique bounded trilinear
extension
\[
T:
L^{p_1}(w_1)\times L^{p_2}(w_2)\times L^{p_3}(w_3)
\longrightarrow L^r(v).
\]
\end{cor}

\begin{rem}
It is natural to ask for the full range of exponents
for which \eqref{eq:bound} holds. A gap remains between
the range established in Theorem~\ref{thm:main} and
the range ruled out by Demeter's counterexample
\cite{Demeter}. The approach in this note extends
the result of \cite{OpenAI}, but determining the
optimal range appears to require new input.
\end{rem}

The rest of the paper is organized as follows.
In Section~\ref{Sec02}, we establish the local estimates
and the corresponding dyadic sparse bound.
In Section~\ref{Sec03}, we pass to the continuous form
and prove Theorem~\ref{prop:continuous}.
Section~\ref{Sec04} contains the proofs of
Theorem~\ref{thm:main} and Corollary~\ref{cor:infty}.
Finally, in Section~\ref{Sec05}, we prove the weighted
estimates in Corollary~\ref{cor:weighted}.

Throughout this note, constants may depend on the
exponents and on the fixed kernel parameters, and,
in the weighted estimates, on the indicated weight
characteristic, but not on the input functions or the scales.

\bigskip 

\section{Local estimates} \label{Sec02}

\subsection{The local input}

We recall the dyadic framework and the kernel class used in
\cite[Section~2]{OpenAI}.

\begin{defn}\label{def:dyadic-lattice}
A \emph{dyadic lattice} on $\R$ is a collection
\[
\calD=\bigcup_{k\in\Z}\calD_k
\]
of bounded half-open intervals such that, for some fixed
$\ell>0$, the following conditions hold:
\begin{enumerate}
\item[(i)]
For every $k\in\Z$, the collection $\calD_k$ is a
partition of $\R$ into intervals of length $\ell 2^{-k}$.

\item[(ii)]
Every interval in $\calD_k$ is the disjoint union of
two intervals in $\calD_{k+1}$.
\end{enumerate}
The two intervals in {\rm (ii)} are called the children
of the original interval. If $I,J\in\calD$ and
$I\subseteq J$, we call $J$ an ancestor of $I$.
\end{defn}

Fix a dyadic lattice $\calD$ and constants $c_0>0$
and $C_0\ge1$.

\begin{defn}\label{def:admissible-kernel}
For $I\in\calD$, a function
$w_I\in C_c^\infty(\R^2)$ is called an \emph{admissible kernel}
if it satisfies the following conditions:
\begin{enumerate}
\item[(i)]
For every $(x,t)\in\supp \; w_I$ and $0\le j\le3$,
\[
x+jt\in I,
\qquad \textrm{and} \qquad 
\dist(x+jt,\partial I)\ge c_0|I|.
\]

\item[(ii)]
For every $u\in\R$ and $0\le j\le3$,
\begin{equation}\label{eq:marginals}
\int_{\R}w_I(u-jt,t)\,dt=0.
\end{equation}

\item[(iii)]
For every multi-index $\alpha$ and every
$(x,t)\in\R^2$,
\[
|\partial^\alpha w_I(x,t)|
\le
|I|^{-|\alpha|}
\bigl(C_0(|\alpha|+2)\bigr)^{C_0(|\alpha|+2)}.
\]
\end{enumerate}
\end{defn}

\begin{rem}
For an explicit construction of admissible kernels, see
the proof of \cite[Lemma~2.4]{OpenAI}.
\end{rem}

For each admissible kernel $w_I$, define the associated
four-linear form by
\begin{equation}\label{eq:HI}
H_I(f_0,f_1,f_2,f_3)
:=
|I|^{-2}
\iint_{\R^2}
w_I(x,t)\prod_{j=0}^3f_j(x+jt)\,dx\,dt.
\end{equation}

\begin{defn}\label{def:rooted-subtree}
A \emph{finite rooted subtree} of $\calD$ with root $J$
is a finite collection $\mathcal T\subset\calD$
containing $J$ such that every $I\in\mathcal T$
is contained in $J$, and
\[
I\in\mathcal T,\quad K\in\calD,\quad
I\subseteq K\subseteq J
\quad\Longrightarrow\quad
K\in\mathcal T.
\]
\end{defn}

\begin{prop}\label{prop:input}
There exists $2<q_*<3$ such that, for every
$q\in(2,q_*)$, the following holds.
Let $\mathcal T$ be a finite rooted subtree
with root $J$.
If
\begin{equation}\label{20261027eq01}
\langle f_j\rangle_{q,I}\le1,
\qquad 0\le j\le3,\quad I\in\mathcal T,
\end{equation}
then, for every subset $\mathcal I\subseteq\mathcal T$,
\begin{equation}\label{eq:normalized}
\sum_{I\in\mathcal I}
|I|\,|H_I(f_0,f_1,f_2,f_3)|
\le B_q|J|.
\end{equation}
The constant $B_q$ depends only on $q$
and the admissibility parameters $c_0$ and $C_0$
in Definition~\ref{def:admissible-kernel}.
In particular, it is independent of the dyadic
lattice, $\mathcal T$, $\mathcal I$, and the
choice of admissible kernels.
\end{prop}

\begin{proof}
The proof combines the normalized estimate in
\cite[Proposition~2.3]{OpenAI} with the stopping
argument in the first part of
\cite[Lemma~2.4]{OpenAI}. Fix $2<q<3$.
Our goal is to bound
\[
\frac1{|J|}
\sum_{I\in\mathcal I}
|I|\,|H_I(f_0,f_1,f_2,f_3)|
\]
uniformly over all finite rooted subtrees
$\mathcal T$ with root $J$, all subsets
$\mathcal I\subseteq\mathcal T$, and all functions
satisfying \eqref{20261027eq01}.
By Definition~\ref{def:admissible-kernel},
{\rm (i)}, we may replace each $f_j$ by
$f_j\mathbf 1_J$ in each form $H_I$,
$I\in\mathcal T$.

For $N\ge1$, let $B_N(q)$ be the least constant
in \eqref{eq:normalized} when $\mathcal T$
is restricted to at most $N$ consecutive generations.
That is, $B_N(q)$ is the least constant such that
\begin{equation}\label{eq:BN-definition}
\sum_{\substack{I\in\mathcal I\\
|I|\ge2^{-(N-1)}|J|}}
|I|\,|H_I(f_0,f_1,f_2,f_3)|
\le B_N(q)|J|
\end{equation}
holds for all finite rooted subtrees $\mathcal T$
with root $J$, all subsets
$\mathcal I\subseteq\mathcal T$, and all functions
satisfying \eqref{20261027eq01}, uniformly over
all dyadic lattices and all families of admissible
kernels with the fixed parameters $c_0$ and $C_0$.

We first observe that $B_N(q)<\infty$.
By Definition~\ref{def:admissible-kernel},
{\rm (i)} and {\rm (iii)}, the Cauchy--Schwarz
inequality, and the changes of variables
\[
(x,t)\mapsto(x,x+t)
\quad\text{and}\quad
(x,t)\mapsto(x+2t,x+3t),
\]
we have
\begin{equation}\label{eq:elementary-local}
|H_I(f_0,f_1,f_2,f_3)|
\lesssim
\prod_{j=0}^3\langle f_j\rangle_{2,I}
\le
\prod_{j=0}^3\langle f_j\rangle_{q,I}
\le1.
\end{equation}
The last inequality follows from
\eqref{20261027eq01}.
The intervals in each fixed generation are
pairwise disjoint subintervals of $J$.
Since the length restriction below permits
at most $N$ generations, we obtain
\[
\sum_{\substack{I\in\mathcal I\\
|I|\ge2^{-(N-1)}|J|}}
|I|\,|H_I(f_0,f_1,f_2,f_3)|
\lesssim N|J|.
\]
Thus $B_N(q)\lesssim N$.

It remains to bound $B_N(q)$ uniformly in $N$.
We shall establish
\[
B_N(q)\le C+\frac12B_N(q),
\]
where $C$ is independent of $N$.
Since $B_N(q)<\infty$, this will imply the desired
uniform bound.
We divide the remainder of the proof into four steps.

\medskip
\noindent
\textbf{Step 1. A local estimate for arbitrary inputs.}

\medskip

\noindent
Fix $N\ge1$.
We first prove the estimate
\begin{equation}\label{eq:local-from-BN}
\sum_{I\in\mathcal A}
|I|\,|H_I(g_0,g_1,g_2,g_3)|
\le B_N(q)|Q|\prod_{j=0}^3a_j,
\end{equation}
where $\mathcal U$ is a finite rooted subtree
with root $Q$ satisfying
\[
|I|\ge2^{-(N-1)}|Q|,
\qquad I\in\mathcal U,
\]
$\mathcal A\subseteq\mathcal U$, and
$a_0,a_1,a_2,a_3>0$ satisfy
\begin{equation}\label{eq:local-normalization}
\langle g_j\rangle_{q,I}\le a_j,
\qquad 0\le j\le3,\quad I\in\mathcal U.
\end{equation}

Indeed, set
\[
h_j:=\frac{g_j}{a_j},
\qquad 0\le j\le3.
\]
By \eqref{eq:local-normalization}, we have
\[
\langle h_j\rangle_{q,I}
=\frac{\langle g_j\rangle_{q,I}}{a_j}
\le1,
\qquad 0\le j\le3,\quad I\in\mathcal U.
\]
By \eqref{eq:BN-definition} and the four-linearity
of $H_I$, we obtain
\begin{align*}
\sum_{I\in\mathcal A}
|I|\,|H_I(g_0,g_1,g_2,g_3)|
&=
\left(\prod_{j=0}^3a_j\right)
\sum_{I\in\mathcal A}
|I|\,|H_I(h_0,h_1,h_2,h_3)|\\
&\le
B_N(q)|Q|\prod_{j=0}^3a_j.
\end{align*}
This proves \eqref{eq:local-from-BN}.

\medskip
\noindent
\textbf{Step 2. A global maximal-function estimate.}

\medskip

\noindent
We next prove that
\begin{equation}\label{eq:aux-global}
\sum_{I\in\mathcal F}
|I|\,|H_I(g_0,g_1,g_2,g_3)|
\lesssim
B_N(q)\int_{\R}\prod_{j=0}^3M_qg_j(x)\,dx
\end{equation}
for all compactly supported $L^q$ functions
$g_0,g_1,g_2,g_3$ and every nonempty finite
collection $\mathcal F\subset\calD$ satisfying
\begin{equation}\label{eq:finite-scale-condition}
\min_{I\in\mathcal F}|I|
\ge
2^{-(N-1)}\max_{I\in\mathcal F}|I|.
\end{equation}
Here
\[
M_qg(x):=\sup_{L\ni x}\langle g\rangle_{q,L},
\]
where the supremum is over bounded intervals $L$.
The implicit constant in \eqref{eq:aux-global}
depends only on $q$, $c_0$, and $C_0$.

Let $\mathcal Q$ be the collection of dyadic
intervals of length $\max_{I\in\mathcal F}|I|$
that contain at least one interval of $\mathcal F$.
The intervals in $\mathcal Q$ are pairwise disjoint,
and every $I\in\mathcal F$ is contained in a unique
member of $\mathcal Q$.

For every dyadic interval $Q$ with
\[
\min_{I\in\mathcal F}|I|
\le |Q|\le
\max_{I\in\mathcal F}|I|,
\]
let $\mathcal C(Q)$ be the collection of maximal
strict dyadic subintervals $Q'\subsetneq Q$
such that
\[
|Q'|\ge\min_{I\in\mathcal F}|I|
\]
and
\begin{equation}\label{eq:stopping-condition}
\langle g_j\rangle_{q,Q'}^q
>
16\langle g_j\rangle_{q,Q}^q
\end{equation}
for at least one $j\in\{0,1,2,3\}$.
By maximality, the intervals in $\mathcal C(Q)$
are pairwise disjoint.
Define
\begin{equation}\label{eq:subtree-Q}
\mathcal T_Q:=
\left\{
I\in\calD:
\begin{array}{l}
I\subseteq Q,\quad
|I|\ge\displaystyle\min_{L\in\mathcal F}|L|,\\[3pt]
I\not\subseteq Q'
\text{ for every }Q'\in\mathcal C(Q)
\end{array}
\right\}.
\end{equation}
The collection $\mathcal T_Q$ contains $Q$
and is closed under passage to dyadic ancestors
inside $Q$.
Hence $\mathcal T_Q$ is a finite rooted subtree
with root $Q$.

By the maximality of the intervals in
$\mathcal C(Q)$, we have
\begin{equation}\label{eq:subtree-averages}
\langle g_j\rangle_{q,I}
\le16^{1/q}\langle g_j\rangle_{q,Q},
\qquad 0\le j\le3,\quad I\in\mathcal T_Q.
\end{equation}
Indeed, any dyadic subinterval $I\subsetneq Q$
violating \eqref{eq:subtree-averages} and satisfying
$|I|\ge\min_{L\in\mathcal F}|L|$ is contained in
a member of $\mathcal C(Q)$, contrary to
\eqref{eq:subtree-Q}.
Moreover, by \eqref{eq:finite-scale-condition},
we have
\[
|I|
\ge\min_{L\in\mathcal F}|L|
\ge2^{-(N-1)}\max_{L\in\mathcal F}|L|
\ge2^{-(N-1)}|Q|,
\qquad I\in\mathcal T_Q.
\]

If $\langle g_j\rangle_{q,Q}>0$ for every $j$,
then by \eqref{eq:local-from-BN}, with
\[
\mathcal U=\mathcal T_Q,\qquad
\mathcal A=\mathcal F\cap\mathcal T_Q,\qquad
a_j=16^{1/q}\langle g_j\rangle_{q,Q},
\]
we have
\begin{equation}\label{eq:subtree-bound}
\sum_{I\in\mathcal F\cap\mathcal T_Q}
|I|\,|H_I(g_0,g_1,g_2,g_3)|
\lesssim
B_N(q)|Q|\prod_{j=0}^3\langle g_j\rangle_{q,Q}.
\end{equation}
If $\langle g_j\rangle_{q,Q}=0$ for some $j$,
then $g_j=0$ almost everywhere on $Q$.
By Definition~\ref{def:admissible-kernel},
{\rm (i)}, we have
\[
H_I(g_0,g_1,g_2,g_3)=0,
\qquad I\in\calD,\quad I\subseteq Q.
\]
Thus \eqref{eq:subtree-bound} also holds when
one of the averages over $Q$ is zero.

To sum \eqref{eq:subtree-bound} over the roots,
define
\begin{equation}\label{eq:sparse-recursion}
\mathcal S_0:=\mathcal Q,\qquad
\mathcal S_{n+1}:=
\bigcup_{Q\in\mathcal S_n}\mathcal C(Q),\qquad
\calS:=\bigcup_{n\ge0}\mathcal S_n.
\end{equation}
Every $Q'\in\mathcal C(Q)$ satisfies
\[
\min_{I\in\mathcal F}|I|
\le |Q'|\le\frac12|Q|.
\]
By \eqref{eq:finite-scale-condition}, we have
$\mathcal S_n=\varnothing$ for $n\ge N$.
In particular, $\calS$ is finite.

For each $I\in\mathcal F$, there is a unique
smallest member $Q\in\calS$ containing $I$.
By \eqref{eq:subtree-Q} and
\eqref{eq:sparse-recursion},
$I\in\mathcal T_Q$ precisely for this choice of $Q$.
Therefore,
\begin{equation}\label{eq:partition-F}
\mathcal F
=
\bigsqcup_{Q\in\calS}
(\mathcal F\cap\mathcal T_Q).
\end{equation}
By \eqref{eq:subtree-bound} and
\eqref{eq:partition-F}, we have
\begin{align}
\sum_{I\in\mathcal F}
|I|\,|H_I(g_0,g_1,g_2,g_3)|
&=
\sum_{Q\in\calS}
\sum_{I\in\mathcal F\cap\mathcal T_Q}
|I|\,|H_I(g_0,g_1,g_2,g_3)|
\nonumber\\
&\lesssim
B_N(q)\sum_{Q\in\calS}|Q|
\prod_{j=0}^3\langle g_j\rangle_{q,Q}.
\label{eq:sparse-sum-bound}
\end{align}

In view of \eqref{eq:sparse-sum-bound},
it remains to prove that
\begin{equation}\label{eq:sparse-maximal-bound}
\sum_{Q\in\calS}|Q|
\prod_{j=0}^3\langle g_j\rangle_{q,Q}
\lesssim
\int_{\R}\prod_{j=0}^3M_qg_j(x)\,dx.
\end{equation}
For $Q\in\calS$ and $0\le j\le3$, define
\begin{equation}\label{eq:stopping-subcollections}
\mathcal C_j(Q):=
\left\{
Q'\in\mathcal C(Q):
\langle g_j\rangle_{q,Q'}^q
>
16\langle g_j\rangle_{q,Q}^q
\right\}.
\end{equation}
If $\langle g_j\rangle_{q,Q}>0$, then, by
the pairwise disjointness of the intervals
in $\mathcal C_j(Q)$, we have
\begin{align*}
\sum_{Q'\in\mathcal C_j(Q)}|Q'|
&\le
\frac{1}{16\langle g_j\rangle_{q,Q}^q}
\sum_{Q'\in\mathcal C_j(Q)}
\int_{Q'}|g_j(x)|^q\,dx\\
&\le
\frac{1}{16\langle g_j\rangle_{q,Q}^q}
\int_Q|g_j(x)|^q\,dx
=\frac{|Q|}{16}.
\end{align*}
If $\langle g_j\rangle_{q,Q}=0$, then
$\langle g_j\rangle_{q,Q'}=0$ for every
$Q'\subseteq Q$, so
$\mathcal C_j(Q)=\varnothing$.
Consequently,
\begin{equation}\label{eq:stopping-measure-j}
\sum_{Q'\in\mathcal C_j(Q)}|Q'|
\le\frac{|Q|}{16},
\qquad 0\le j\le3.
\end{equation}
Since
$\mathcal C(Q)=\bigcup_{j=0}^3\mathcal C_j(Q)$,
we obtain
\begin{equation}\label{eq:stopping-measure}
\sum_{Q'\in\mathcal C(Q)}|Q'|
\le
\sum_{j=0}^3
\sum_{Q'\in\mathcal C_j(Q)}|Q'|
\le\frac14|Q|.
\end{equation}

For $Q\in\calS$, set
\begin{equation}\label{eq:sparse-remainders}
E_Q:=Q\setminus
\bigcup_{Q'\in\mathcal C(Q)}Q'.
\end{equation}
By \eqref{eq:stopping-measure}, we have
\[
|E_Q|\ge\frac34|Q|,
\qquad Q\in\calS.
\]
The sets $E_Q$, $Q\in\calS$, are pairwise disjoint.
Indeed, if $R,Q\in\calS$ and $R\subsetneq Q$,
then, by \eqref{eq:sparse-recursion},
$R\subseteq Q'$ for some $Q'\in\mathcal C(Q)$.
Hence $E_R\cap E_Q=\varnothing$.
If $R$ and $Q$ are not nested, then
$R\cap Q=\varnothing$.

For $Q\in\calS$, we have
\[
M_qg_j(x)\ge\langle g_j\rangle_{q,Q},
\qquad x\in E_Q,\quad 0\le j\le3.
\]
It follows that
\[
|Q|\prod_{j=0}^3\langle g_j\rangle_{q,Q}
\lesssim
|E_Q|\prod_{j=0}^3\langle g_j\rangle_{q,Q}
\le
\int_{E_Q}\prod_{j=0}^3M_qg_j(x)\,dx.
\]
By the pairwise disjointness of the sets $E_Q$,
we obtain
\[
\sum_{Q\in\calS}|Q|
\prod_{j=0}^3\langle g_j\rangle_{q,Q}
\lesssim
\sum_{Q\in\calS}
\int_{E_Q}\prod_{j=0}^3M_qg_j(x)\,dx
\le
\int_{\R}\prod_{j=0}^3M_qg_j(x)\,dx.
\]
This proves \eqref{eq:sparse-maximal-bound}.
By \eqref{eq:sparse-sum-bound} and
\eqref{eq:sparse-maximal-bound}, we obtain
\eqref{eq:aux-global}, uniformly in $N$,
the dyadic lattice, and the admissible kernels.

\medskip
\noindent
\textbf{Step 3. A normalized estimate with a small
coefficient.}

\medskip

\noindent
We next prove that there exists $2<q_*<3$ such that,
for every $q\in(2,q_*)$ and every sufficiently small
$\varepsilon>0$,
\begin{equation}\label{eq:imported-normalized}
\sum_{\substack{I\in\mathcal I\\
|I|\ge2^{-(N-1)}|J|}}
|I|\,|H_I(f_0,f_1,f_2,f_3)|
\le (C+\varepsilon K)|J|,
\end{equation}
where $\mathcal T$ is a finite rooted subtree
with root $J$, $\mathcal I\subseteq\mathcal T$,
the functions $f_0,f_1,f_2,f_3$ satisfy
\eqref{20261027eq01}, and $K$ satisfies
\begin{equation}\label{eq:global-hypothesis}
\sum_{I\in\mathcal F}
|I|\,|H_I(g_0,g_1,g_2,g_3)|
\le
K\int_{\R}\prod_{j=0}^3M_qg_j(x)\,dx
\end{equation}
for every nonempty finite collection $\mathcal F$
satisfying \eqref{eq:finite-scale-condition}
and all compactly supported $L^q$ functions
$g_0,g_1,g_2,g_3$, uniformly over the dyadic
lattices and admissible kernels.
The constant $C$ in \eqref{eq:imported-normalized}
depends only on $q$, $\varepsilon$, $c_0$, and $C_0$;
in particular, it is independent of $N$ and $K$.

We first specify the range of exponents permitted
by the proof of \cite[Proposition~2.3]{OpenAI}.
The restriction on $q$ in
\cite[Section~10.3, p.~92]{OpenAI} has the form
\begin{equation}\label{eq:q-smallness}
\gamma(q-2)<\frac18,
\end{equation}
where $\gamma>0$ is independent of $q$.
Indeed, the discussion immediately following
the parameter table in
\cite[Section~10.3, p.~92]{OpenAI}
specifies that the complexity exponents entering
$\gamma$ are fixed before $q$ is chosen.
Choose $2<q_*<3$ such that
\[
\gamma(q_*-2)<\frac18.
\]
Then
\[
\gamma(q-2)<\gamma(q_*-2)<\frac18,
\qquad 2<q<q_*.
\]
For each fixed $q\in(2,q_*)$, the remaining
auxiliary parameters in
\cite[Sections~10.3--10.4]{OpenAI}
are chosen after $q$.
Thus the proof of \cite[Proposition~2.3]{OpenAI}
is valid for every $q\in(2,q_*)$.

Fix such a $q$ and define
\[
\mathcal T_N:=
\left\{
I\in\mathcal T:
|I|\ge2^{-(N-1)}|J|
\right\}.
\]
The collection $\mathcal T_N$ is a finite rooted
subtree with root $J$ and at most $N$ generations.
Moreover, by \eqref{20261027eq01}, we have
\[
\langle f_j\rangle_{q,I}\le1,
\qquad 0\le j\le3,\quad I\in\mathcal T_N.
\]

The least constant in \eqref{eq:global-hypothesis}
is at most $K$.
By \cite[Proposition~2.3, pp.~4--5]{OpenAI},
with the tree $\mathcal T_N$ and the subset
$\mathcal I\cap\mathcal T_N$, we have
\begin{align*}
\sum_{\substack{I\in\mathcal I\\
|I|\ge2^{-(N-1)}|J|}}
|I|\,|H_I(f_0,f_1,f_2,f_3)|
&=
\sum_{I\in\mathcal I\cap\mathcal T_N}
|I|\,|H_I(f_0,f_1,f_2,f_3)|\\
&\le (C+\varepsilon K)|J|.
\end{align*}
This proves \eqref{eq:imported-normalized}
with the stated dependence of $C$.

\medskip
\noindent
\textbf{Step 4. A uniform bound for $B_N(q)$.}

\medskip

\noindent
Fix $q\in(2,q_*)$.
We finally prove that
\begin{equation}\label{eq:BN-uniform}
B_N(q)\le2C,
\qquad N\ge1,
\end{equation}
where $C$ depends only on $q$, $c_0$, and $C_0$.

By \eqref{eq:aux-global}, for each $N$ we may
choose a constant $K$ in
\eqref{eq:global-hypothesis} satisfying
\[
K\lesssim B_N(q).
\]
The implicit constant depends only on
$q$, $c_0$, and $C_0$, and is independent of
$N$ and $\varepsilon$.
We may therefore fix $\varepsilon>0$ sufficiently
small, depending only on $q$, $c_0$, and $C_0$,
such that
\begin{equation}\label{eq:absorption-choice}
\varepsilon K\le\frac12B_N(q)
\end{equation}
for every $N$ and the corresponding choice of $K$.

By \eqref{eq:imported-normalized} and
\eqref{eq:absorption-choice}, we have
\[
\frac1{|J|}
\sum_{\substack{I\in\mathcal I\\
|I|\ge2^{-(N-1)}|J|}}
|I|\,|H_I(f_0,f_1,f_2,f_3)|
\le C+\varepsilon K
\le C+\frac12B_N(q).
\]
Taking the supremum over all normalized data
in \eqref{eq:BN-definition}, we obtain
\begin{equation}\label{eq:BN-absorption}
B_N(q)\le C+\frac12B_N(q).
\end{equation}
Since $B_N(q)<\infty$, it follows that
\[
\frac12B_N(q)\le C,
\]
which proves \eqref{eq:BN-uniform}.
The choice of $\varepsilon$ depends only on
$q$, $c_0$, and $C_0$, so $C$ has precisely
the dependence asserted above.

Finally, let $\mathcal T$ be a finite rooted
subtree with root $J$.
Choose $N$ large enough that
\[
|I|\ge2^{-(N-1)}|J|,
\qquad I\in\mathcal T.
\]
By \eqref{eq:BN-definition} and
\eqref{eq:BN-uniform}, we have
\[
\sum_{I\in\mathcal I}
|I|\,|H_I(f_0,f_1,f_2,f_3)|
\le B_N(q)|J|
\le2C|J|
\]
for every $\mathcal I\subseteq\mathcal T$
and all functions satisfying \eqref{20261027eq01}.
This proves \eqref{eq:normalized} with $B_q=2C$.
\end{proof}

\begin{rem}
OSVerify\footnote{See
\url{https://osverify.github.io/cases/trilinear-hilbert.html},
Issue~01.}
points out that the nesting assumption in
\cite[Lemma~8.6]{OpenAI} omits the inclusion
$P_{N,h}\subseteq C_{N,h}$ for $0\le h\le J$.
This inclusion holds automatically in the application in
\cite[Proposition~8.9]{OpenAI}, since the PRE graph
refines the CURRENT graph; see \cite[Section~8.4]{OpenAI}.
Adding this inclusion to the hypotheses makes the
counting argument valid without changing the chain-length
bound or the subsequent normalized estimate in
\cite[Proposition~2.3]{OpenAI} used above.
\end{rem}

\subsection{The stopping time argument}

By Proposition~\ref{prop:input} and
\eqref{eq:BN-definition}, we have
\[
B_N(q)\le B_q,
\qquad N\ge1,\quad 2<q<q_*.
\]
We use this uniform bound in
\eqref{eq:sparse-sum-bound}, together with the sparse
family constructed in Step~2 of the proof of
Proposition~\ref{prop:input}, to obtain the following local sparse domination estimate.

\begin{prop}\label{prop:local-sparse}
Fix $2<q<q_*$.
For every finite collection $\mathcal I\subset\calD$,
every family of admissible kernels, and all
$f_0,f_1,f_2,f_3\in L^q_{\mathrm{loc}}(\R)$,
there exists a finite $3/4$-sparse collection
$\calS\subset\calD$ such that
\begin{equation}\label{eq:local-sparse}
\sum_{I\in\mathcal I}
|I|\,|H_I(f_0,f_1,f_2,f_3)|
\le
C_q\Lambda_{q,\calS}(f_0,f_1,f_2,f_3).
\end{equation}
The constant $C_q$ depends only on $q$, $c_0$,
and $C_0$.
\end{prop}

\begin{proof}
The case $\mathcal I=\varnothing$ is immediate.
Otherwise, choose $N\ge1$ so that
\eqref{eq:finite-scale-condition} holds with
$\mathcal F=\mathcal I$.

With $\mathcal F=\mathcal I$ and $g_j=f_j$,
let $\calS$ be the finite $3/4$-sparse family
constructed in Step~2 of the proof of
Proposition~\ref{prop:input}.
The decomposition \eqref{eq:sparse-recursion}
and the estimate \eqref{eq:sparse-sum-bound}
involve only the restrictions of the functions
to the finitely many dyadic intervals of length
$\max_{I\in\mathcal I}|I|$ that contain at least
one interval of $\mathcal I$.
Hence they remain valid for
$f_j\in L^q_{\mathrm{loc}}(\R)$.

By \eqref{eq:sparse-sum-bound} and
$B_N(q)\le B_q$, we have
\begin{align*}
\sum_{I\in\mathcal I}
|I|\,|H_I(f_0,f_1,f_2,f_3)|
&\lesssim
B_N(q)\Lambda_{q,\calS}(f_0,f_1,f_2,f_3)\\
&\le
B_q\Lambda_{q,\calS}(f_0,f_1,f_2,f_3).
\end{align*}
This proves \eqref{eq:local-sparse}.
\end{proof}

\bigskip 

\section{The continuous form and a sparse domination estimate} \label{Sec03}

In this section, we prove Theorem~\ref{prop:continuous}.
The proof combines the sparse estimate in
Proposition~\ref{prop:local-sparse} with the
kernel construction in the proof of
\cite[Lemma~2.4]{OpenAI}.

In what follows, we choose the fixed admissibility parameters
$c_0$ sufficiently small and $C_0$ sufficiently
large.
Moreover, let $q_*$ be the exponent in Proposition~\ref{prop:input} corresponding to
the fixed admissibility parameters $c_0$ and $C_0$.

\begin{proof}[Proof of Theorem~\ref{prop:continuous}]
\medskip
\noindent
\textbf{Step 1. A representation of the continuous kernel.}

\medskip 

\noindent
We first establish the identity
\begin{equation}\label{eq:kernel-representation}
\int_0^\infty
\sum_{m\in\Z}W\left(\frac{x}{a}-m, \frac{t}{a} \right)\,\frac{da}{a^2}
=
-\frac1t,
\qquad x\in\R,\quad t\ne0,
\end{equation}
where $W\in C_c^\infty(\R^2)$ has four vanishing
line marginals.

Choose an odd compactly supported Gevrey function
$\varphi$ and a compactly supported Gevrey function
$\omega$ such that\footnote{Here, a compactly supported
smooth function $g$ is called Gevrey if there exist
constants $A,B>0$ and $s>1$ such that
$\|g^{(n)}\|_{L^\infty(\R)}\le AB^n(n!)^s$
for every integer $n\ge0$. The parameter $s$ is called
the Gevrey order.}
\[
\varphi'''(0)=1,
\qquad
\sum_{m\in\Z}\omega(X-m)=1.
\]
Define
\begin{equation}\label{eq:W-definition}
W(X,S):=
\prod_{j=0}^3(\partial_S-j\partial_X)
[\omega(X)\varphi(S)].
\end{equation}
For $0\le j\le3$, set
\[
W_j(X,S):=
\prod_{\substack{0\le k\le3\\k\ne j}}
(\partial_S-k\partial_X)
[\omega(X)\varphi(S)].
\]
Then
\[
W(u-jS,S)
=
\frac{d}{dS}W_j(u-jS,S).
\]
Since $W_j$ is compactly supported, we have
\begin{equation}\label{eq:W-marginals}
\int_{\R}W(u-jS,S)\,dS=0,
\qquad u\in\R,\quad 0\le j\le3.
\end{equation}

By \eqref{eq:W-definition} and the identity
$\sum_m\omega(X-m)=1$, we have
\begin{equation}\label{eq:spatialsum}
\sum_{m\in\Z}W\left(\frac{x}{a}-m, \frac{t}{a} \right)
=
\varphi^{(4)} \left(\frac{t}{a} \right),
\qquad a>0.
\end{equation}
For $t>0$, the substitution $s=t/a$ yields
\[
\int_0^\infty
\varphi^{(4)} \left(\frac{t}{a} \right)\,\frac{da}{a^2}
=\frac1t\int_0^\infty\varphi^{(4)}(s)\,ds
=-\frac{\varphi'''(0)}t
=-\frac1t.
\]
Since $\varphi^{(4)}$ is odd, we obtain
\begin{equation}\label{eq:recover}
\int_0^\infty
\varphi^{(4)} \left(\frac{t}{a} \right)\,\frac{da}{a^2}
=
-\frac1t,
\qquad t\ne0.
\end{equation}
By \eqref{eq:spatialsum} and \eqref{eq:recover},
we have \eqref{eq:kernel-representation}.

\medskip
\noindent
\textbf{Step 2. A uniform estimate for finitely many scales.}

\medskip 

\noindent
We next prove that
\begin{equation}\label{eq:continuous-finite-scales}
|\Lambda_{M,b}(f_0,f_1,f_2,f_3)|
\lesssim
C_q\sup_{\calS}
\Lambda_{q,\calS}(f_0,f_1,f_2,f_3),
\qquad M\ge1,\quad 1\le b<2,
\end{equation}
where
\begin{equation}\label{eq:Lambda-Mb}
\begin{split}
\Lambda_{M,b}(f_0,f_1,f_2,f_3)
:={}&
\sum_{n=-M}^{M}\frac1{b2^n}
\iint_{\R^2}
\varphi^{(4)}\left(\frac{t}{b2^n}\right)
\prod_{j=0}^3f_j(x+jt)\,dx\,dt.
\end{split}
\end{equation}
The implicit constant will be independent of $M$
and $b$.

Choose $R\ge1$ such that
\[
\supp \; W\subset[-R,R]^2,
\]
and fix a dyadic number $D\ge 100R$.
For $1\le b<2$ and
$\theta\in\{0,1/3,2/3\}$, define
\[
\calD_b^\theta
:=
\bigcup_{n\in\Z}
\left\{
Db2^n
\bigl[k+(-1)^n\theta,k+1+(-1)^n\theta\bigr):
k\in\Z
\right\}.
\]
Each $\calD_b^\theta$ is a dyadic lattice.
The alternating signs ensure nesting because
$3\theta\in\Z$.

Fix $a=b2^n$ and $m\in\Z$.
On the support of $W\left(\frac{x}{a}-m, \frac{t}{a} \right)$, we have
\[
x+jt\in a[m-4R,m+4R],
\qquad 0\le j\le3.
\]
At scale $Da$, observe that the boundaries of the three
lattices have offsets separated by $Da/3$.
We may therefore choose an interval
$I=I_{n,m,b}$ in one of the lattices
$\calD_b^\theta$ such that
\[
|I|=Da,
\qquad
am\in I,
\qquad
\dist(am,\partial I)\ge\frac{Da}{3}.
\]
Consequently,
\[
x+jt\in I,
\qquad
\dist(x+jt,\partial I)
\ge\frac{Da}{3}-4Ra
\ge\frac{|I|}{4},
\qquad 0\le j\le3,
\]
on the support of $W\left(\frac{x}{a}-m, \frac{t}{a} \right)$.

For fixed $n$, the points $am$, $m\in\Z$,
are separated by $a$.
Thus each interval $I$ receives at most $D+1$
values of $m$.
We may partition the pairs $(n,m)$ into at most
$3(D+1)$ classes such that, within each class,
the assigned intervals are distinct and belong
to a single dyadic lattice.

Within each class, define
\begin{equation}\label{eq:continuous-admissible}
w_I(x,t):=
D\,W\left(\frac{x}{a}-m, \frac{t}{a} \right),
\qquad I=I_{n,m,b},\quad a=b2^n.
\end{equation}
The support condition in
Definition~\ref{def:admissible-kernel}, {\rm (i)},
holds for $c_0\le1/4$.
By \eqref{eq:W-marginals} and the substitution
$t=aS$, we have
\[
\int_{\R}w_I(u-jt,t)\,dt
=
Da\int_{\R}W(u/a-m-jS,S)\,dS
=
0.
\]
Moreover, since Gevrey classes are closed under products
and differentiation, the fixed function $W$ satisfies
\[
\|\partial^\alpha W\|_\infty
\le C^{|\alpha|+1}(|\alpha|!)^s
\]
for some fixed $C\ge1$ and $s>1$.
Since $|I|=Da$, we obtain
\[
\begin{aligned}
|\partial^\alpha w_I(x,t)|
&\le
|I|^{-|\alpha|}(CD)^{|\alpha|+1}(|\alpha|!)^s\\
&\le
|I|^{-|\alpha|}
\bigl(C_0(|\alpha|+2)\bigr)^{C_0(|\alpha|+2)}.
\end{aligned}
\]
The last inequality follows from
$|\alpha|!\le(|\alpha|+2)^{|\alpha|}$
by choosing $C_0$ sufficiently large,
depending only on the fixed constants $C,D,s$.
Thus the kernels in
\eqref{eq:continuous-admissible} belong to a
single admissible class, independently of
$n$, $m$, and $b$.

By \eqref{eq:HI} and
\eqref{eq:continuous-admissible}, we have
\begin{equation}\label{eq:continuous-atom}
\frac1a
\iint_{\R^2}
W\left(\frac{x}{a}-m, \frac{t}{a} \right)\prod_{j=0}^3f_j(x+jt)\,dx\,dt
=
|I|H_I(f_0,f_1,f_2,f_3).
\end{equation}
For fixed $M$ and $b$, only finitely many
spatial translates contribute, since $f_0$
is compactly supported.
By Proposition~\ref{prop:local-sparse} in each
of the finitely many classes, we obtain
\[
\sum_{n=-M}^{M}\sum_{m\in\Z}
\left|
\frac1{b2^n}
\iint_{\R^2}
W\left(\frac{x}{b2^n}-m,\frac{t}{b2^n}\right)
\prod_{j=0}^3f_j(x+jt)\,dx\,dt
\right|
\lesssim
C_q\sup_{\calS}
\Lambda_{q,\calS}(f_0,f_1,f_2,f_3).
\]
By \eqref{eq:spatialsum} and
\eqref{eq:Lambda-Mb}, this proves
\eqref{eq:continuous-finite-scales}.

\medskip
\noindent
\textbf{Step 3. Passage to the principal-value form.}

\medskip 

\noindent
We finally prove that
\begin{equation}\label{eq:continuous-limit}
\lim_{M\to\infty}
\int_1^2
\Lambda_{M,b}(f_0,f_1,f_2,f_3)\,\frac{db}{b}
=
\Lambda(f_0,f_1,f_2,f_3).
\end{equation}
By \eqref{eq:continuous-finite-scales}, we have
\begin{equation}\label{eq:continuous-integrated-bound}
\left|
\int_1^2
\Lambda_{M,b}(f_0,f_1,f_2,f_3)\,\frac{db}{b}
\right|
\lesssim
C_q\sup_{\calS}
\Lambda_{q,\calS}(f_0,f_1,f_2,f_3),
\end{equation}
uniformly in $M$.

Set
\[
K_M(t):=
\int_{2^{-M}}^{2^{M+1}}
\varphi^{(4)} \left(\frac{t}{a} \right)\,\frac{da}{a^2}.
\]
The substitution $a=b2^n$ in
\eqref{eq:Lambda-Mb} gives
\begin{equation}\label{eq:integrated-kernel}
\int_1^2
\Lambda_{M,b}(f_0,f_1,f_2,f_3)\,\frac{db}{b}
=
\iint_{\R^2}
K_M(t)\prod_{j=0}^3f_j(x+jt)\,dx\,dt.
\end{equation}
Since $\varphi^{(4)}$ is odd, $K_M$ is odd.
Furthermore,
\begin{equation}\label{eq:kernel-limit-bound}
K_M(t)\longrightarrow-\frac1t,
\qquad
|K_M(t)|
\le
\int_0^\infty
|\varphi^{(4)} \left(\frac{t}{a} \right)|\,\frac{da}{a^2}
\lesssim\frac1{|t|},
\qquad t\ne0.
\end{equation}
The convergence follows from \eqref{eq:recover}.

Write
\[
G(x,t):=\prod_{j=1}^3f_j(x+jt).
\]
By the oddness of $K_M$, the right-hand side of
\eqref{eq:integrated-kernel} equals
\[
\int_{\R}f_0(x)
\int_0^\infty
K_M(t)\bigl(G(x,t)-G(x,-t)\bigr)\,dt\,dx.
\]
For $x\in\supp f_0$, the smoothness and decay
of the Schwartz functions imply
\[
|G(x,t)-G(x,-t)|
\lesssim
\begin{cases}
t, & 0<t\le1,\\
t^{-2}, & t>1,
\end{cases}
\]
uniformly in $x$.
Together with \eqref{eq:kernel-limit-bound},
this provides an integrable majorant independent
of $M$.
By dominated convergence, we obtain
\begin{align*}
\lim_{M\to\infty}
\int_1^2
\Lambda_{M,b}(f_0,f_1,f_2,f_3)\,\frac{db}{b}
&=
-\int_{\R}f_0(x)
\int_0^\infty
\frac{G(x,t)-G(x,-t)}{t}\,dt\,dx\\
&=
\Lambda(f_0,f_1,f_2,f_3).
\end{align*}
This proves \eqref{eq:continuous-limit}.

Finally, combining \eqref{eq:continuous-integrated-bound}
and \eqref{eq:continuous-limit}, we obtain
\eqref{eq:continuous}.
\end{proof}

\medskip 

\section{Proof of the $L^p$ estimates} \label{Sec04}

We now deduce Theorem~\ref{thm:main} and
Corollary~\ref{cor:infty} from the sparse bound in
Theorem~\ref{prop:continuous}.

\begin{proof}[Proof of Theorem~\ref{thm:main}]
Fix $p_1,p_2,p_3$ and $r$ as in Theorem~\ref{thm:main}.
Choose $q\in(2,q_*)$ sufficiently close to $2$ that
\begin{equation}\label{eq:chooseq}
q<\min\{p_1,p_2,p_3\},
\qquad
r<q':=\frac{q}{q-1}.
\end{equation}
Such a choice is possible because $p_j>2$, $r<2$,
and $q'\to2$ as $q\downarrow2$.
We first consider $f_1,f_2,f_3\in C_c^\infty(\R)$.

We first treat the Banach range $1\le r<2$.
Recall that
\[
M_qf=\bigl(M(|f|^q)\bigr)^{1/q},
\]
where $M$ is the uncentered Hardy--Littlewood maximal
operator.
For a $3/4$-sparse collection $\calS$, let
$E_I\subset I$, $I\in\calS$, be pairwise disjoint
measurable sets satisfying $|E_I|\ge3|I|/4$.
Since $\langle f_j\rangle_{q,I}\le M_qf_j(x)$ for
$x\in E_I$, we have
\begin{align}
\Lambda_{q,\calS}(f_0,f_1,f_2,f_3)
&\lesssim
\sum_{I\in\calS}
\int_{E_I}\prod_{j=0}^3M_qf_j(x)\,dx
\nonumber\\
&\le
\int_{\R}\prod_{j=0}^3M_qf_j(x)\,dx.
\label{eq:maximal}
\end{align}

Set $p_0=r'$ if $r>1$ and $p_0=\infty$ if $r=1$.
By \eqref{eq:chooseq}, we have
\[
p_j>q,
\qquad 0\le j\le3,
\qquad
\sum_{j=0}^3\frac1{p_j}=1.
\]
By \eqref{eq:continuous}, \eqref{eq:maximal},
H\"older's inequality, and the boundedness of $M_q$
on $L^{p_j}$, we obtain
\[
|\Lambda(f_0,f_1,f_2,f_3)|
\lesssim
\prod_{j=0}^3\norm{M_qf_j}_{L^{p_j}}
\lesssim
\prod_{j=0}^3\norm{f_j}_{L^{p_j}}.
\]
Taking the supremum over bounded compactly supported
$f_0$ with $\norm{f_0}_{L^{p_0}}\le1$, we obtain
\eqref{eq:bound} for $1\le r<2$.

Next, we treat the quasi-Banach range $2/3<r<1$.
By $2<q<q_*<3$ and \eqref{eq:chooseq}, we have
\[
q\ge1,\qquad \frac4q>1,\qquad
q<p_j<\infty,\quad 1\le j\le3,\qquad r<q'.
\]
By \cite[Corollary~2.15]{LMO} and
\eqref{eq:continuous}, with three inputs, sparse
exponent vector $(q,q,q,q)$, sparsity constant $3/4$,
and all weights equal to one, we obtain
\eqref{eq:bound} for $2/3<r<1$.

Finally, by a standard density argument, $T$ extends
uniquely to a bounded trilinear map from
$L^{p_1}(\R)\times L^{p_2}(\R)\times L^{p_3}(\R)$
into $L^r(\R)$ satisfying \eqref{eq:bound}.
\end{proof}

\begin{proof}[Proof of Corollary~\ref{cor:infty}]
The three cases are handled in the same way, and here 
we only give the proof for $p_1=\infty$.
Choose $q\in(2,q_*)$ such that
\[
q<\min\{p_2,p_3,r'\}.
\]
Since
\[
\frac1{r'}+\frac1{p_2}+\frac1{p_3}=1,
\qquad
\norm{M_qf_1}_{L^\infty}\le\norm{f_1}_{L^\infty},
\]
by \eqref{eq:continuous}, \eqref{eq:maximal}, and
H\"older's inequality, we have
\begin{equation}\label{eq:inftyform}
|\Lambda(f_0,f_1,f_2,f_3)|
\lesssim
\norm{f_0}_{L^{r'}}
\norm{f_1}_{L^\infty}
\norm{f_2}_{L^{p_2}}
\norm{f_3}_{L^{p_3}}
\end{equation}
for $f_0,f_1,f_2,f_3\in C_c^\infty(\R)$.

Define
\[
{\mathbb A}(f_0,f_2,f_3)(y)
:=
\operatorname{p.v.}\int_{\R}
f_0(y+t)f_2(y-t)f_3(y-2t)\,\frac{dt}{t}.
\]
The smoothness and compact support of the inputs
ensure that the principal value defines a locally
integrable function.
By the change of variables $y=x-t$ in the truncated
integrals and passage to the principal value, we have
\begin{equation}\label{eq:adjoint-pairing}
\Lambda(f_0,f_1,f_2,f_3)
=
\int_{\R}f_1(y){\mathbb A}(f_0,f_2,f_3)(y)\,dy
\end{equation}
for every $f_1\in C_c^\infty(\R)$.
By \eqref{eq:adjoint-pairing} and \eqref{eq:inftyform},
we obtain
\begin{align}
\norm{{\mathbb A}(f_0,f_2,f_3)}_{L^1}
&=
\sup_{\substack{f_1\in C_c^\infty(\R)\\
\norm{f_1}_{L^\infty}\le1}}
|\Lambda(f_0,f_1,f_2,f_3)|
\nonumber\\
&\lesssim
\norm{f_0}_{L^{r'}}
\norm{f_2}_{L^{p_2}}
\norm{f_3}_{L^{p_3}}.
\label{eq:adjoint-L1}
\end{align}
Since $r',p_2,p_3<\infty$, by standard density argument, ${\mathbb A}$ extends
to a bounded trilinear map
\[
{\mathbb A}:L^{r'}(\R)\times L^{p_2}(\R)\times L^{p_3}(\R)
\longrightarrow L^1(\R).
\]

For $f_1\in L^\infty$, $f_2\in L^{p_2}$, and
$f_3\in L^{p_3}$, define $T(f_1,f_2,f_3)\in L^r$ by
\begin{equation}\label{eq:infty-extension}
\int_{\R}f_0(x)T(f_1,f_2,f_3)(x)\,dx
=
\int_{\R}f_1(y){\mathbb A}(f_0,f_2,f_3)(y)\,dy
\end{equation}
for every $f_0\in L^{r'}$.
By \eqref{eq:adjoint-L1}, \eqref{eq:adjoint-pairing},
and duality, \eqref{eq:infty-extension}
defines the required bounded trilinear extension of $T$.
\end{proof}

\section{Proof of Corollary \ref{cor:weighted}} \label{Sec05}

Finally, we prove the weighted estimates for the trilinear Hilbert transform $T$. 

\begin{proof}[Proof of Corollary~\ref{cor:weighted}]
Let $p_1,p_2,p_3,r$ and $\vec w$ be as in
Corollary~\ref{cor:weighted}, and set
\[
K:=[\vec w]_{A_{\vec p,\vec 2}},
\qquad
v:=\prod_{j=1}^3w_j^{r/p_j}.
\]
Define
\[
u_0:=v^{\frac{2}{2-r}},
\qquad
u_j:=w_j^{-\frac{2}{p_j-2}},
\quad 1\le j\le3.
\]
By \cite[Lemma~5.3]{LMO} and  duality, we have
\[
[u_0]_{A_{\frac{2r}{2-r}}}
\le K^{\frac{2r}{2-r}},
\qquad
[u_j]_{A_{\frac{2p_j}{p_j-2}}}
\le K^{\frac{2p_j}{p_j-2}},
\quad 1\le j\le3,
\]
where $A_s$ denotes the classical Muckenhoupt class.
The reverse H\"older inequality therefore provides
$\delta>0$, depending only on the exponents and $K$,
such that
\begin{equation}\label{eq:weighted-rh}
\langle u_i\rangle_{1+\delta,I}
\lesssim
\langle u_i\rangle_{1,I},
\qquad 0\le i\le3,
\end{equation}
for every bounded interval $I\subset\R$, with an
absolute implicit constant.

Choose $q\in(2,q_*)$ sufficiently close to $2$ that
\begin{equation}\label{eq:weighted-exponents}
q<\min\{p_1,p_2,p_3\},
\qquad
r<q':=\frac{q}{q-1},
\end{equation}
and
\begin{equation}\label{eq:weighted-q-choice}
1<\theta_i\le1+\delta,
\qquad 0\le i\le3,
\end{equation}
where
\[
\theta_0
:=
\frac{q'(2-r)}{2(q'-r)},
\qquad
\theta_j
:=
\frac{q(p_j-2)}{2(p_j-q)},
\quad 1\le j\le3.
\]
Such a choice is possible because
$\theta_i\to1$ as $q\downarrow2$ for every
$0\le i\le3$.

Set $\vec q=(q,q,q,q)$.
In the notation of \cite{LMO}, the corresponding
weight characteristic is
\[
[\vec w]_{A_{\vec p,\vec q}}
:=
\sup_I
\bigl\langle v^{\frac{q'}{q'-r}}\bigr\rangle_{1,I}
^{\frac1r-\frac1{q'}}
\prod_{j=1}^3
\bigl\langle w_j^{-\frac{q}{p_j-q}}\bigr\rangle_{1,I}
^{\frac1q-\frac1{p_j}}.
\]
By \eqref{eq:weighted-rh},
\eqref{eq:weighted-q-choice}, and H\"older's inequality,
we have
\begin{align}
[\vec w]_{A_{\vec p,\vec q}}
&=
\sup_I
\langle u_0\rangle_{\theta_0,I}^{\frac1r-\frac12}
\prod_{j=1}^3
\langle u_j\rangle_{\theta_j,I}^{\frac12-\frac1{p_j}}
\nonumber\\
&\lesssim
\sup_I
\langle u_0\rangle_{1,I}^{\frac1r-\frac12}
\prod_{j=1}^3
\langle u_j\rangle_{1,I}^{\frac12-\frac1{p_j}}
=K.
\label{eq:weighted-improvement}
\end{align}
In particular, $\vec w\in A_{\vec p,\vec q}$.

Since \eqref{eq:weighted-exponents} holds and
$4/q>1$, by Theorem~\ref{prop:continuous},
\eqref{eq:weighted-improvement}, and
\cite[Corollary~2.15]{LMO}, we have
\[
\norm{T(f_1,f_2,f_3)}_{L^r(v)}
\le
C\prod_{j=1}^3\norm{f_j}_{L^{p_j}(w_j)}
\]
for $f_1,f_2,f_3\in C_c^\infty(\R)$.
Both $\delta$ and $q$ may be chosen in terms of
the exponents and $K$, so the constant $C$ has
the dependence asserted in
Corollary~\ref{cor:weighted}.
By a standard density argument, $T$ extends
uniquely to the stated weighted spaces.
\end{proof}

\medskip

\noindent{\bf AI disclosure.}
The author used ChatGPT for language polishing and proofreading.
The main analytic input in Section~\ref{Sec02}, which is
the central part of this note, comes from the core argument
of OpenAI's work \cite{OpenAI} (entry 086 in the OpenAI
Research Catalog).
The author proposed investigating whether a sparse
domination estimate could be extracted from that proof
to extend OpenAI's boundedness result.
The adaptation of the stopping-time argument in
Sections~\ref{Sec02} and \ref{Sec03} was developed through repeated
discussions with ChatGPT, which also assisted with
drafting and revising those sections.
The initial draft of Sections~\ref{Sec02} and \ref{Sec03} generated by ChatGPT
was largely unreadable and required substantial reorganization
and repeated revision under the author's direction.
The remaining mathematical contributions to this note
are due to the author.
The author takes full responsibility for the content
and conclusions of this paper.

\medskip 

\noindent{\bf Acknowledgement.}
The author was supported by the NSF grant DMS-2555999
and by the Simons Travel grant MPS-TSM-00007213.

\end{document}